\documentclass[11pt]{article}
\usepackage{amssymb}

\newcommand{\F}{{\mathbb F}}
\newcommand{\EE}{{\mathbb E}}

\newcommand{\LL}{\Lambda}
\newcommand{\eps}{\varepsilon}
\newcommand{\SIGMA}{\raisebox{-0.4ex}{\mbox{\Large $\Sigma$}}}

\newtheorem{theorem}{Theorem}
\newtheorem{lemma}{Lemma}
\newtheorem{proposition}{Proposition}

\title{Subsets of $\F_{p^n}$ without three term arithmetic progressions 
have several large Fourier coefficients}

\author{Ernie Croot}

\begin{document}

\maketitle

\begin{abstract}  Suppose that $f : \F_{p^n} \to [0,1]$ satisfies 
$$
\SIGMA_a f(m)\ =\ \theta F\ \in\ [F^{8/9}, F],\ {\rm where\ }F\ =\ |\F_{p^n}|\ =\ p^n.
$$
In this paper we will show the following:  Let $f_j$ denote the size of the $j$th
largest Fourier coefficient of $f$.  If 
$$
f_j\ <\ \theta^{j^{1/2 + \delta}} F,
$$
for some integer $j$ satisfying
$$
J_0(\delta,p)\ <\ j\ <\ F^{1/8},
$$
then $S = {\rm support}(f)$ contains a non-trivial three-term arithmetic progression.  
Thus, the result is asserting that if the Fourier transform decays rapidly enough
(though not all that rapidly -- in particular, not quite exponentially fast), then 
$S$ is forced to have a three-term arithmetic progression.  This result is similar
in spirit to that appearing in \cite{croot}; however, in that paper the focus was on 
the ``small'' Fourier coefficients, whereas here the focus is on the ``large'' Fourier
coefficients (furthermore, the proof in the present paper requires much more
sophisticated arguments than those of that other paper).

Here is a  partial description of how this result was proved:  First, to get our proof
started, we reduced to the case where $\theta$, $\delta$ and $j$
satisfied certain nice constraints.  For example, Meshulam's theorem \cite{meshulam} 
was used to reduce to the case $\theta < 1/p$; then, we reduced to the case 
$\delta < 2/3$ by using a ``dimension collapsing and cutting argument'' from
\cite{croot} which says that for special ``smooth'' functions $f$, the underlying set
$S$ must always be rich in three-term arithmetic progressions; and finally, this same
argument (dimension collapsing) was used to reduce to the case where $j < n^{2-\delta}$. 
The rest of the proof  used a type of Roth-Meshulam \cite{meshulam}
iteration.  Unfortunately, because our main theorem is to work with 
densities $\theta$ that can be much smaller than $n^{-1}$, 
which is the limit of the basic Roth-Meshulam \cite{meshulam}
approach, we cannot expect to use the usual ``density increment'' principle to
reach a subset of $S$ containing lots of three-term arithmetic progressions.
Instead, we showed that the number 
of ``large Fourier coefficients'' decreases by a lot at each iteration.  
To show that this count indeed decreases 
by ``a lot'', and not just by one (as one can get by a trivial argument), 
we had to unravel some of the additive structure of the set of large Fourier coefficients. 
One type of argument that could perhaps do this for us is that of 
Shkredov \cite{shkredov}; unfortunately, it appears that his argument does not
work well in our present context.  Instead, a ``phase shifting and pigeonhole'' 
argument, which appears in Lemma \ref{prop1}, was used to get at some of this additive 
structure.   But the hypotheses of this lemma are slightly unusual, and in order
to make the Lemma useful for proving our main theorem, we needed to
introduce an extra possibility into our Roth-Meshulam \cite{meshulam} 
iteration scheme; the two
possibilities are basically the last two possible conclusions of Proposition 
\ref{prop0}.  Because these conclusions are somewhat different
from each other, we needed to introduce certain ``invariants'', which are 
(\ref{inv0}) through (\ref{inv6}), in order to show that our process 
eventually terminates with $S$ having a non-trivial three-term arithmetic progression.     
\end{abstract}

We will assume throughout that 
$$
\F\ :=\ \F_{p^n},\ {\rm and\ } F\ :=\ p^n\ =\ |\F|.
$$

Suppose that 
$$
f\ :\ \F\ \to\ [0,1],
$$
and that
$$
\EE(f)\ =\ F^{-1} \SIGMA_a f(a)\ =\ \theta.
$$
Define
$$
\LL(f)\ :=\ \EE_{m,d}(f(m)f(m+d)f(m+2d))\ =\ F^{-2} \SIGMA_{m,d} 
f(m)f(m+d)f(m+2d).
$$
Note that if $f$ is an indicator function for some set, then $\LL(f)$ is some normalized
count of the number of three-term arithmetic progressions in this set.
\bigskip

In this paper we will prove a theorem which says that if $f$ has too few large
Fourier coefficients, then 
$$
\LL(f)\ >\ \theta F^{-1},
$$
which would thus imply that the set $S$ given by
$$
S\ :=\ {\rm support}(f)\ =\ \{ m \in \F\ :\ f(m)\ >\ 0\},
$$
contains three-term arithmetic progressions.
\bigskip

In order to properly state this result, we will need a few more definitions:  
First, given a function $f$ defined on $\F$, and given $a \in \F$, we let 
$$
\hat f(a)\ :=\ \SIGMA_m f(m) e^{2\pi i (a \cdot m) / p},
$$
where the $a\cdot m$ is the usual dot product.  Next, write
$$
\F\ =\ \{a_1,...,a_F\},
$$
where the $a_i$s are ordered so that 
$$
|\hat f(a_1)|\ \geq\ |\hat f(a_2)|\ \geq\ \cdots\ \geq\ |\hat f(a_F)|.
$$
Let us set 
$$
a_1\ =\ 0
$$
(there may be other values of $a$ such that $|\hat f(a)| = |\hat f(0)|$), and let us 
assume that
$$
a_2 = -a_3,\ a_4 = -a_5,\ ...,
$$
which we know is possible because from the fact that $f : \F \to [0,1]$ we 
deduce that 
$$
\hat f(a)\ =\ \overline{\hat f(-a)}.
$$
\bigskip

The main theorem of this paper will show that if the Fourier transform of $f$
decays rapidly enough, in the sense that $|\hat f(a_j)|$ is small enough, then
we can deduce that $S$ contains a non-trivial three-term arithmetic progression.
This sort of result was proved in \cite{croot}; however, the focus of that paper was
more on properties of the ``small'' Fourier coefficients, whereas in the present
paper we work only with the ``large'' Fourier coefficients.   Our theorem is 
as follows:

\begin{theorem} \label{main_theorem}  For any $\delta > 0$, and $p$ prime,
the following holds for all dimensions $n$ (of our field $\F_{p^n}$) 
sufficiently large:  Suppose that 
$$
F^{-1/9}\ <\ \theta\ \leq\ 1,
$$
and suppose that for some integer $j$ satisfying
$$
J_0(\delta,p)\ <\ j\ <\ F^{1/8},
$$
we have that 
\begin{equation} \label{faj}
|\hat f(a_j)|\ <\ \theta^{j^{1/2+\delta}} F.
\end{equation}
Then, $S$ contains a three-term arithmetic 
progression; more specifically, we have that
$$
\LL(f)\ >\ \theta F^{-1}.
$$
\end{theorem}
\bigskip

\noindent {\bf Remark 1.}  We note that a key strength of this theorem is 
this exponent $1/2 + \delta$ in (\ref{faj}); for, if the exponent were replaced 
by the larger value $1 + o(1)$, then the proof would be profoundly easier,
and indeed Proposition \ref{prop3} below gives such a result.   

One way to get a feel for what this result is saying is to suppose that 
$\theta \sim n^{-100}$ (much too small to be dealt with using Roth-Meshulam 
\cite{meshulam} iteration), and then suppose that $f$ has at most, say, $n^{1.9}$ (the $1.9$ can
be any number smaller than $2$) Fourier coefficients $\hat f(a)$ satisfying
$$
F^{2/3}\ \leq\ |\hat f(a)|\ \leq\ \theta F.
$$
Thus, a dyadic interval 
$$
[2^{-t} \theta F,\ 2^{-t+1} \theta F]\ \subseteq\ [F^{2/3},\ \theta F]
$$
contains on average about $n^{0.9 - o(1)}$ Fourier coefficients $|\hat f(a)|$ -- or rather,
norms of Fourier coefficients.  Note that our Theorem \ref{main_theorem} tells
us that in this case $S$ contains a non-trivial three-term arithmetic progression.
\bigskip

\noindent {\bf Remark 2.}  We feel that Theorem \ref{main_theorem} should have a
much simpler proof than we give here in the present paper, perhaps along the
lines of the argument in \cite{croot}; furthermore, we 
believe that it ought to be possible to prove a much, much stronger result.  Here are 
two conjectures along these lines, the second much stronger than the first:
\bigskip

\noindent {\bf Conjecture 1.}  There exists $0 < c < 1$ such that the following holds
for every $\delta > 0$, $p \geq 3$ prime, and dimensions $n$ sufficiently large:
If 
$$
\theta\ >\ F^{-c},\ {\rm and\ } |\hat f(a_j)| < \theta^{j^\delta} F,
$$
for some
$$
J_0(\delta,p)\ <\ j\ <\ F^{1/8},
$$ 
then $S$ contains a non-trivial three-term arithmetic progression.  
\bigskip

\noindent {\bf Conjecture 2.}  There exists $0 < c_1 < 1$ and $c_2 > 0$ 
such that the following holds for all $p \geq 3$ prime, 
and dimensions $n$ sufficiently large:  If
$$
\theta\ >\ F^{-c_1},\ {\rm and\ } |\hat f(a_j)| < \theta^{c_2 \log j} F,
$$
for some
$$
J_0(p)\ <\ j\ <\ F^{1/8},
$$
then $S$ contains a non-trivial three-term arithmetic progression.
\bigskip

\noindent {\bf Remark 3.}  Here is a  partial description of how the proof of 
Theorem \ref{main_theorem} goes:  First, to get the proof started we 
reduce to some case where $\theta, \delta$ and $j$ satisfy certain nice 
constraints.  For example, we use Meshulam's theorem \cite{meshulam} 
to reduce to the case
$\theta < 1/p$; then, we reduce to the case $\delta < 2/3$ by using a 
``dimension collapsing and cutting argument'' from \cite{croot} which says that
for special ``smooth'' functions $f$, the underlying set $S$ must always be 
rich in three-term arithmetic progressions; and, finally, we use this same 
argument (dimension collapsing) to reduce to the case $j < n^{2-\delta}$.  
The rest of the proof uses a type of Roth-Meshulam \cite{meshulam}
iteration.  Unfortunately, because Theorem \ref{main_theorem} is to work 
with densities $\theta$ that can
be much smaller than $n^{-1}$, which is the limit of the basic Roth-Meshulam 
\cite{meshulam} approach, we cannot expect to use the usual ``density increment'' principle to
reach a subset of $S$ containing lots of three-term arithmetic progressions.
To get around this, we instead show that at each Roth-Meshulam \cite{meshulam}
iteration the number of ``large Fourier coefficients'' decreases by a lot.  To show that this count indeed 
decreases by ``a lot'', and not just by one (as one can get by a trivial argument), 
we must unravel some of the additive structure of the set of large Fourier coefficients. 
One type of argument that could perhaps do this for us is that of 
Shkredov \cite{shkredov}; unfortunately, it appears that his argument does not
work well in our present context.  Instead, we use a ``phase shifting and pigeonhole'' 
argument, which appears in Lemma \ref{prop1} below, to get at some of this additive 
structure.   But the hypotheses of this Lemma are slightly unusual, and in order
to make the Lemma useful for proving Theorem \ref{main_theorem}, we need to
introduce an extra possibility in our Roth-Meshulam \cite{meshulam} 
iteration scheme; the two
possibilities are basically the last two possible conclusions of Proposition 
\ref{prop0}, which is stated in the next section.  
Because these conclusions are somewhat different
from each other, we need to introduce certain ``invariants'', which are (\ref{inv0})
through (\ref{inv6}), in order to show that our process eventually terminates with
$S$ having a non-trivial three-term arithmetic progression.     
\bigskip

\section{Proof of Theorem \ref{main_theorem}}

\subsection{Preliminary results}

The proof of our theorem will follow a certain type of Roth iteration, but one which
incorporates some of the non-trivial additive structure of sets having no
three-term arithmetic progressions.  The main proposition we use which 
incorporates this additive structure into Roth iteration is as follows.

\begin{proposition} \label{prop0}
For every prime $p \geq 3$, there exist integers
$$
n_0(p)\ {\rm and\ }\ j_0(p)\ >\ 1
$$
such that if the following all hold
\bigskip

$\bullet$ 
$$
n\ =\ {\rm dim}(\F)\ >\ n_0(p),\ j\ >\ j_0(p);
$$
\bigskip

$\bullet$
$$
f\ :\ \F\ \to\ [0,1],\ {\rm where\ } \EE(f)\ =\ \theta\ >\ 2 F^{-1/8};\ {\rm and,}
$$
\bigskip

$\bullet$
\begin{equation} \label{faj_bound}
|\hat f(a_j)|\ <\ \theta^{j^{1/2 + \delta}} F,\ {\rm and\ } |\hat f(a_{j-1})|\ >\ 2 F^{-1/2},
\end{equation}
\bigskip

\noindent then one of the following must hold:
\bigskip

$\bullet$  Either $S = {\rm support}(f)$ contains a non-trivial three-term arithmetic
progression;
\bigskip

$\bullet$ or, for 
$$
\ell\ =\ \lceil j/50 \rceil
$$
we have that 
$$
|\hat f(a_{\ell})|\ <\ \theta^{2\ell^{1/2 + \delta}} F;
$$
\bigskip

$\bullet$ or, there exists a function $h\ :\ \F_{p^{n-1}}\ \to\ [0,1]$ such that 
$$
{\rm support}(h)\ {\rm has\ 3APs}\ \ \Longrightarrow\ \ S\ {\rm has\ 3APs},
$$
such that
$$
\EE(h)\ =\ p^{-n+1} \SIGMA_{w \in \F_{p^{n-1}}} h(w)\ \geq\ \theta,
$$
and such that if we let
$$
\{b_1,...,b_{p^{n-1}}\}\ =\ \F_{p^{n-1}}\ \ {\rm satisfy\ \ }
|\hat h(b_1)|\ \geq\ \cdots\ \geq\ |\hat h(b_{p^{n-1}})|,
$$
and set
$$
b_1\ =\ 0,\ b_2 = -b_3,\ b_4 = -b_5,\ ...,
$$
then for
$$
\ell\ =\ j\ +\ 51\ -\ \lfloor (j/50)^{1/2} \rfloor
$$
we will have
$$
|\hat h(b_{\ell})|\ \leq\ \theta^{j^{1/2 + \delta}} p^n.
$$
\end{proposition}

Next, we will need the following two Propositions to clean up certain ``exceptional
cases'' that arise later in the body of the proof of Theorem \ref{main_theorem}.
Furthermore, the second of these, Proposition \ref{prop3}, 
can be thought of as a weak version of 
Theorem \ref{main_theorem}, where instead of having $j^{1/2}$ in
the exponent $\theta^{j^{1/2 + \delta}}$, we have $j + 2$.

\begin{proposition} \label{prop30} Suppose $f : \F \to [0,1]$, that our dimension
$n$ is sufficiently large in terms of $p$, that 
$$
\EE(f)\ >\ 2pF^{-1/8},
$$
and that for some 
$$
2\ \leq\ j\ <\ F^{1/8}
$$
we have 
$$
|\hat f(a_j)|\ <\ \theta^2/32.
$$
Then, $S = {\rm support}(f)$ contains a non-trivial three-term arithmetic progression.
\end{proposition}

\begin{proposition} \label{prop3}  Suppose $f: \F \to [0,1]$, that
$$
n\ =\ {\rm dim}(\F)\ \geq\ n_0(p),
$$
where $n_0(p)$ is some function only of the characteristic $p$ of the field $\F$,
and that
$$
\EE(f)\ =\ \theta\ >\ 2p F^{-1/8}.
$$
If for some 
$$
2\ \leq\ j\ \leq\ F^{1/8}
$$
we have that 
$$
|\hat f(a_j)|\ <\ \theta^{j+2} F/2,
$$
then $S = {\rm support}(f)$ contains a non-trivial three-term arithmetic progression.
\end{proposition}

\subsection{Body of the Proof of Theorem \ref{main_theorem}}

\subsubsection{Initial reductions}

First, we can assume that 
$$
\theta\ <\ \min(1/p,1/4),
$$ 
for if $\theta \geq 1/p$ or $1/4$, then by Meshulam's theorem \cite{meshulam} we will 
have that once $n$ is sufficiently large
(as a functions of $p$ alone), the set $S$ contains a three-term arithmetic
progression, thereby proving Theorem \ref{main_theorem}.

Next, we may assume that 
\begin{equation} \label{hit_me}
\theta^2/32\ <\ |\hat f(a_j)|\ <\ \theta^{j^{1/2 + \delta}} F,
\end{equation}
since otherwise failure of this lower bound to hold would imply, 
by Proposition \ref{prop30}, that $S$ contains a
three-term arithmetic progression, which again would prove  
Theorem \ref{main_theorem}.

From this it follows that we may assume that 
\begin{equation} \label{theta_lower}
\theta\ >\ F^{-2j^{-1/2}}\ >\ 2pF^{-1/8},
\end{equation}
(with $j > 256$) for if this first inequality does not hold, then 
we will have from the hypotheses of Theorem 
\ref{main_theorem} that 
$$
|\hat f(a_j)|\ <\ \theta^{j^{1/2 + \delta}} F\ \leq\ F^{-2 j^\delta} F\ \leq\ F^{-1}.
$$
So, since $\theta > F^{-1/8}$, it follows that this is smaller than $\theta^2/32$
(for $n$ sufficiently large), which would thus imply that $S$ contains a three-term
arithmetic progression.

\subsubsection{Moving to a benign index $j'$}

Now suppose that we are unlucky and have that  
$$
|\hat f(a_{j-1})|\ <\ 2F^{1/2}.
$$
If this holds, then let $j' < j$ be the largest index for which 
$$
|\hat f(a_{j'-1})|\ \geq\ 2F^{1/2},\ {\rm but\ } |\hat f(a_{j'})|\ <\ 2F^{1/2}.
$$
Such an index $j'$ clearly exists, since 
$$
\hat f(0)\ =\ \theta F\ >\ 2F^{7/8}\ >\ 2 F^{1/2}.
$$

Next we claim that, in view of (\ref{theta_lower}), 
we may assume that $j'$ is as large as we might happen to
need, simply by choosing $J_0$ as large as we like, where $J_0$
is assumed to be a lower bound for $j$; more specifically, we will
have that if 
$$
j'\ \leq\ J_0^{1/2}/6\ \leq\ j^{1/2}/6,
$$ 
then $S$ would have to contain a three-term arithmetic progression:  
To see this, note that
$$
|\hat f(a_{j'})|\ <\ 2 F^{1/2}\ <\ (F^{-2j^{-1/2}})^{j' + 2} F/2\ <\ \theta^{j'+2} F/2.
$$
Thus, Proposition \ref{prop3} would imply that for $n$ sufficiently large, the 
set $S$ contains three-term arithmetic progressions.  
\bigskip

Next, observe that (\ref{hit_me}) implies that for $j$ sufficiently large,
$$
|\hat f(a_{j'})|\ <\ 2F^{1/2}\ <\ \theta^{j^{1/2 + \delta}/4}F\ <\ \theta^{j^{1/2 + \delta/2}} F
\ <\ \theta^{(j')^{1/2 + \delta/2}} F.
$$

What this means is that we have passed from the pair
$$
(j,\delta)\ \to\ (j',\delta'),\ {\rm with\ } \delta' = \delta/2,
$$
such that for this new instance of $j$ and $\delta$ we have the hypotheses of
Theorem \ref{main_theorem} hold, at least if we choose $J_0$ large enough
so that $j' > j_0$, because we have that
\begin{equation} \label{this_here}
\theta^{(j')^{1/2 + \delta'}} F\ >\ 2 F^{1/2},\ {\rm and\ }|\hat f(a_{j'-1})|\ \geq\ 2 F^{1/2}.
\end{equation}
The reason that this is useful is that it will allow us to use Proposition 
\ref{prop0}; futhermore, we will apply Proposition \ref{prop0} iteratively, and
at each step of the iteration we will want that the invariants (\ref{this_here}) are
maintained (as well as several other invariants listed below).
\bigskip

For notational convenience, we will assume that we have (\ref{this_here}) holding
with $j'$ replaced with $j$, and $\delta'$ replaced with $\delta$.  In other words,
we will assume
\begin{equation} \label{this_here2}
\theta^{ j^{1/2 + \delta}} F\ >\ 2F^{1/2},\ {\rm and\ } |\hat f(a_{j-1})|\ \geq\ 2 F^{1/2}.
\end{equation}

\subsubsection{Further reductions}

We may assume that 
$$
\delta\ \leq\ 2/3,
$$
at least for $j \geq 7$, since if $\delta > 2/3$, we would have that
$$
|\hat f(a_j)|\ <\ \theta^{j^{1/2 + \delta}} F\ <\ \theta^{j^{7/6}} F/2,
$$
and then Proposition \ref{prop3} would imply that $S$ contains a three-term
arithmetic progression.

Furthermore, we may assume that 
$$
j\ \leq\ n^{2 -\delta},
$$
for if $j > n^{2 - \delta}$, then using the facts that $\delta < 2/3$ and that 
$\theta < 1/p$, along with the hypotheses of our theorem, we would have 
that  
$$
|\hat f(a_j)|\ <\ \theta^{j^{(1/2 + \delta)(2 - \delta)}} F\ =\ p^{-n^{1 + 3\delta/2 - \delta^2}}F 
<\ \theta^2/32,
$$
which by our arument above near (\ref{hit_me}) would imply that, again, $S$
contains a three-term arithmetic progression.

\subsubsection{Invariants of applying Proposition \ref{prop0} iteratively}

We now apply Proposition \ref{prop0} iteratively, amplifying the value of 
$\delta > 0$ at each step, until we get ourselves in a position where we 
can apply Proposition \ref{prop3}.  We will think of each such iteration as a process
that takes a particular instance of $j$ and $\delta$, and produces a new instance; so,
application of Proposition \ref{prop0} iteratively corresponds to a sequence
$$
(j_1,\delta_1)\ :=\ (j,\delta)\ \to\ (j_2,\delta_2)\ \to\ (j_3,\delta_3)\ \to\ \cdots
$$
At each iteration, rather than having a function $f : \F \to [0,1]$, we will have 
a function $h_i : \F_{p^n_i} \to [0,1]$; in other words, we also will have a corresponding
sequence of functions and dimensions given by
$$
(h_1,n_1)\ :=\ (f,n)\ \to\ (h_2,n_2)\ \to\ (h_3,n_3)\ \to \cdots,
$$
where 
$$
n_{i+1}\ =\ n_i\ {\rm or\ } n_i-1.
$$
We get that $n_{i+1} = n_i-1$ precisely if the last bullet of the conclusion of
Proposition \ref{prop0} is applied; if the next-to-last bullet is appled, we instead get
$n_{i+1} = n_i$, and just set $h_{i+1} = h_i$.
\bigskip

The process continues until we reach $(j_T, \delta_T)$ satisfying
\begin{equation} \label{the_end}
{\rm either\ }j_T\ <\ j_0,\ {\rm or\ } \delta_T\ >\ 2/3,
\end{equation}
whichever of these occurs first, 
where $j_0 = j_0(p)$ is as appears in Proposition \ref{prop0}.
\bigskip

At each process of the iteration we need to maintain a number of invariants, in 
order for Proposition \ref{prop0} to apply at the next iteration, and in order for
us to later show that we reach a $\delta_T > 2/3$.  These invariants are:
\bigskip

$\bullet$ We have that if $\{b_1,...,b_{p^{n_i}}\}$ are the elements of 
$\F_{p^{n_i}}$ arranged so that 
$$
|\hat h_i(b_1)|\ \geq\ |\hat h_i(b_2)|\ \geq\ \cdots,
$$
with
$$
b_1 = 0,\ b_2 = -b_3,\ b_4  = -b_5, ...,
$$
then
\begin{equation} \label{inv0}
|\hat h_i(b_{j_i})|\ <\ \theta^{j_i^{1/2 + \delta_i}} p^{n_i}.
\end{equation}
\bigskip

$\bullet$  The ineqalities in (\ref{this_here2}) will be maintained; that is,
\begin{equation} \label{inv1}
\theta^{ j_i^{1/2 + \delta_i}} p^{n_i}\ >\ 2p^{n_i/2},\ {\rm and\ } 
|\hat f(a_{j_i-1})|\ \geq\ 2 p^{n_i/2}.
\end{equation}
\bigskip

$\bullet$ We have that
\begin{equation} \label{inv2}
\EE(h_i)\ =\ p^{-n_i} \SIGMA_a h_i(a)\ \geq\ \theta.
\end{equation}
\bigskip

$\bullet$ We have that
\begin{equation} \label{inv3}
{\rm support}(h_i)\ {\rm has\ a\ 3AP}\ \Longrightarrow\ S\ {\rm has\ a\ 3AP}.
\end{equation}
\bigskip

$\bullet$ We will have that, except possibly at the last iteration, all the
\begin{equation} \label{inv4}
j_i > j_0.
\end{equation}  
\bigskip

$\bullet$ We have that, except at the last iteration, 
\begin{equation} \label{inv41}
\delta_i\ \leq\ 2/3.
\end{equation}
\bigskip

$\bullet$ The dimensions $n_i$ never get too small; in fact, we will have that 
\begin{equation} \label{inv5}
n_i > 99n/100.
\end{equation}
\bigskip

$\bullet$ A final invariant, which we will prove in a later subsection is that
\begin{equation} \label{inv6}
{1/2 + \delta_{i+1} \over 1/2 + \delta_i}\ \geq\ 1 + {\log\log(j_i) - \log\log(j_{i+1}) \over 10}.
\end{equation}
\bigskip

From this last invariant, it will follows that 
$$
{1/2 + \delta_T \over 1/2 + \delta}\ \geq\ 1 + {\log\log(j) - \log\log(j_T) \over 10}.
$$
And so, if our process terminates with $j_T < j_0$ we will either have that
$$
{1/2 + \delta_T \over 1/2 + \delta}\ \geq\ 1 + {\log\log(J_0) - \log\log(j_0) \over 10},
$$
which can be made  as large as desired by choosing $J_0$ as large as we 
like (relative to $j_0$); or else, our process terminates early with $\delta_T > 2/3$,
but with $j > j_0$.  Either way, we can have our process end with 
$$
\delta_T \geq 3\ {\rm and\ } j_T > j_0,\ {\rm or\ with\ } \delta_T > 2/3\ {\rm but\ } 
2 \leq j_T \leq j_0. 
$$
Either way, if $j_0$ is large enough, it will give us that
$$
|\hat h_T(b_{j_T})|\ <\ \theta^{j_T + 2} p^{n_T}/2,
$$
and therefore since by (\ref{inv5}) we have $n_T > n/2$, it follows also that
$$
\EE(h_T)\ \geq\ \theta\ >\ p^{-n/9}\ >\ 2 p^{-n_T/8},
$$   
Proposition \ref{prop3} will imply then that the support of $h_T$
contains a three-term progression, which would mean that $S$ does as well.

\subsubsection{Proving that the invariants all hold}

Let us suppose that we have already applied Proposition \ref{prop0} $t-1$ times,
and so have produced
$$
(j_1, \delta_1)\ \to\ (j_2,\delta_2)\ \to\ \cdots\ \to\ (j_t,\delta_t),
$$
as well as the corresponding functions $h_i$ and the dimensions $n_i$.
We will assume that 
$$
j_t \geq j_0(p),
$$
since otherwise $t = T$ and we are done.  
\bigskip

We note now that the hypotheses of Proposition \ref{prop0} hold for the
function $f$ replaced with $h_t$ (and $\F$ by $\F_{p^{n_i}}$).  So, one or the
other of the conclusions of that proposition must hold.  We may assume, moreover,
that one of the last two conclusions must holds, since otherwise $S$ has 
three-term progressions and we are done.
\bigskip

\noindent {\bf Case 1 (next to last conclusion).}
Suppose that the next to last conclusion of Proposition \ref{prop0} holds.  Then, we
set 
\begin{equation} \label{h_equal}
h_{t+1}\ :=\ h_t,\  n_{t+1}\ :=\ n_t,\ j_{t+1}\ :=\ \lceil j_t/50 \rceil,
\end{equation}
and we let $\delta_{t+1}$ satisfy
\begin{equation} \label{h_equal2}
2 j_{t+1}^{1/2 + \delta_t}\ =\ j_{t+1}^{1/2 + \delta_{t+1}}.
\end{equation}
Note that our conclusion of Proposition \ref{prop0} gives us
$$
|\hat h_{t+1}(b_{t+1})|\ <\ \theta^{j_{t+1}^{1/2 + \delta_{t+1}}} p^{n_{t+1}},
$$
which means that (\ref{inv0}) holds.

Since $h_t$ satisfies the second part of invariant (\ref{inv1}), and since
we have (\ref{h_equal}), we must have that
\begin{eqnarray*}
|\hat h_{t+1}(b_{j_{t+1}-1})|\ \geq\ |\hat h_{t+1}(b_{j_{t+1}})|\ =\ 
|\hat h_t(b_{\lceil j_t/50 \rceil})| 
&\geq&\ |\hat h_t(b_{j_t-1})| \\
&\geq&\ 2 p^{n_{t+1}/2}.
\end{eqnarray*}
So, the second part of (\ref{inv1}) holds for $h_{t+1}$.   
Furthermore,
$$
2 p^{n_{t+1}/2}\ \leq\ |\hat h_{t+1}(b_{j_{t+1}})|\ <\ 
\theta^{j_{t+1}^{1/2 + \delta_{t+1}}} p^{n_{t+1}}
$$
implies that the first part of (\ref{inv1}) holds for $h_{t+1}$; so, we have that 
(\ref{inv1}) holds for our new function $h_{t+1}$ and
choice of $j_{t+1}$ and $\delta_{t+1}$.  

From the fact that $h_{t+1} = h_t$ we get for free that 
(\ref{inv2}), (\ref{inv3}), and (\ref{inv5}) all hold.  
The only invariant we have to show, then, to finish this case is that 
(\ref{inv6}) holds, which we now do (we don't have to worry about 
(\ref{inv4}) or (\ref{inv41}) 
because the only time they could be violated is at the last iteration):  
We note from (\ref{h_equal2}) that
$$
{1/2 + \delta_{t+1} \over 1/2 + \delta_t}\ =\ 1 + {\log 2 \over (1/2 + \delta_t)\log j_{t+1}}.
$$
On the other hand, we also have that 
\begin{eqnarray*}
{\log\log(j_t) - \log\log(j_{t+1}) \over 10}\ &\leq&\ {\log\log(j_t) - \log \left ( (\log j_t) 
\left ( 1 - {\log 50 \over \log j_t} \right ) \right ) \over 10}\\
&=&\ - {\log (1 - \log(50)/\log j_t) \over 10} \\
&\leq&\ {\log 50 \over 5 \log j_t} \\
&\leq&\ {\log 50 \over 5 \log j_{t+1}} \\
&\leq&\ {\log 2 \over (1/2 + \delta_t) \log j_{t+1}},
\end{eqnarray*}
at least for $j_t$ sufficiently large and $\delta_t \leq 2/3$, both of which we assume to
be true.  So, we have that the invariant
(\ref{inv6}) holds.
\bigskip

\noindent {\bf Case 2 (the last conclusion).}   
Suppose that the last conclusion of Proposition \ref{prop0} holds.  Then, we
set 
\begin{equation} \label{case2_info}
h_{t+1}\ :=\ h,\ n_{t+1}\ :=\ n_t-1,\ {\rm and\ } j_{t+1}\ =\ j_t + 51 - \lfloor (j_t/50)^{1/2} 
\rfloor,
\end{equation}
where $h$ is as given in the proposition.  One of the conclusions of 
Proposition \ref{prop0} can thus be restated as
$$
|\hat h_{t+1}(b_{j_{t+1}})|\ <\ \theta^{j_t^{1/2 + \delta_t}} p^{n_t}.
$$
We let $\delta_{t+1} > 0$ be defined by
\begin{equation} \label{equal_theta}
\theta^{j_{t+1}^{1/2 + \delta_{t+1}}} p^{n_{t+1}}\ =\ \theta^{j_t^{1/2 + \delta_t}}
p^{n_t}.
\end{equation}
The fact that such $\delta_{t+1} > 0$ exists and satisfies
$$
\delta_{t+1}\ >\ \delta_t,
$$
is guaranteed by the facts that $\theta < 1/p$, $n_{t+1} = n_t-1$, and 
and $j_t > j_0$.

Now we have that 
$$
|\hat h_{t+1}(b_{j_{t+1}})|\ <\ \theta^{j_{t+1}^{1/2 + \delta_{t+1}}} p^{n_{t+1}},
$$
which therefore means that (\ref{inv0}) holds.

From the fact that $h_t$ satisfied (\ref{inv1}) we have that
$$
\theta^{j_t^{1/2 + \delta_t}} p^{n_t}\ \geq\ 2 p^{n_t/2}\ >\ 2 p^{n_{t+1}/2},
$$
which, along with (\ref{equal_theta}), implies that  
$$
|\hat h_{t+1}(b_{j_{t+1}-1})|\ >\ 2 p^{n_{t+1}/2}.
$$
So, we have that (\ref{inv1}) holds as well for $h_{t+1}$.

Furthermore, (\ref{inv2}), (\ref{inv3}), and (\ref{inv41}) will all hold.  We will
have to hold off for the time being on showing that (\ref{inv4}) holds, as its
proof amounts to showing that it does not take more than $n/100$ iterations before
our process of constructing the functions $h_i$ terminates.

It remains to show that (\ref{inv6}) holds.  To do this, we observe from
(\ref{equal_theta}) that 
\begin{equation} \label{bigo}
{1/2 + \delta_{t+1} \over 1/2 + \delta_t}\ =\ {\log j_t \over \log j_{t+1}}\ +\ 
O \left ( {\log p \over j_t^{1/2+\delta_t} \log j_{t+1}} \right ).
\end{equation}
Now we claim that 
$$
{\log j_t \over \log j_{t+1}}\ >\ 1 + \log\log j_t - \log\log j_{t+1},
$$
which can be seen by letting $j_{t+1} = j_t^{1-\gamma}$, $0 < \gamma < 1$, 
and then noting that this inequality is equivalent to
$$
{1 \over 1 - \gamma}\ >\ 1 - \log(1 - \gamma),
$$
which is easy to verify on taking a Taylor expansion.

So, to verify (\ref{inv6}) we just need to address this big-oh error term above.
This we do by noting that the value of $j_{t+1}$ given in 
(\ref{case2_info})  implies that 
$$
{\log j_t \over \log j_{t+1}}\ >\ {1 \over 1 + {1 \over \log j_t} 
\log \left ( 1 - { (j_t/50)^{1/2} - 51 \over j_t}  \right )}.
$$
By having $j_0$ sufficiently large, we can arrage to have the right-hand-side
exceed 
$$
1\ +\ {1 \over 8 j_t^{1/2} \log j_t};
$$
and, in fact, by choosing $j_0$ large enough, we can have that the 
big-oh error term on (\ref{bigo}) will be strictly smaller than 
$$
{1 \over 2} \left ( {\log j_t \over \log j_{t+1}} - 1 \right ).
$$
So, it follows that 
$$
{1/2 + \delta_{t+1} \over 1/2 + \delta_t}\ >\ 1\ +\ {\log\log j_t - \log\log j_{t+1} \over 2},
$$
which thus establishes (\ref{inv6}) for $h_{t+1}$ (in fact, it establishes somewhat
more).

We note that three still remains the problem of showing that (\ref{inv5}) holds,
and we will establish this in the next sub-sub-section.

\subsubsection{A lower bound on the residual dimension, and the conclusion of
the proof}

To fix this last loose end of showing that (\ref{inv5}) holds, note that
only case 2 above, where $n_{t+1} = n_t - 1$, could cause us problems.   
The absolute worst thing that could happen, then, is if we were in Case 2 
every single step of the way.   Let us see what value for
the final dimension $n_T$ this would give:  First we claim that the absolute
most number of times we could pass through Case 2 is:
\begin{equation} \label{m_bound}
T\ <\ 20 j^{1/2} \log j,
\end{equation}
which would prove that at each iteration,
$$
n_t\ >\ n - 20 j^{1/2} \log j\ >\ n - n^{1 - \delta/2 + o(1)}\ >\ 99 n/100.
$$ 
Here, we have used the fact that 
$$
n\ <\ j^{2 - \delta}.
$$
\bigskip

Let us now see that (\ref{m_bound}) holds:  First, we will use the fact that
$$
j_{t+1}\ <\ j_t\ -\ j_t^{1/2}/10
$$ 
at least if $j_t > j_0$ is sufficiently large.  So, if we run the 
process Case 2 for at least $10j^{1/2}$ steps, then we claim that 
we will reach a 
$$
j_t\ <\ j/2.
$$
To see this, note that at each step $i$ where  
$$
i\ <\ t\ :=\ \lfloor 10 j^{1/2} \rfloor,
$$ 
if we always had that 
$$
j_i\ >\ j/2,
$$
then we would get that at step $t$ that 
\begin{equation} \label{jtdown}
j_t\ <\ j\ -\ (t-1)(j/2)^{1/2}/10\ <\ j/2,
\end{equation}
for $j$ sufficiently large (say $j > j_0$).  

Applying (\ref{jtdown}) iteratively, we see that so long as $j'$ is sufficiently large,
if we pass through Case 2 for 
$$
m\ =\ \lfloor  10 j^{1/2} \log(j/j')/\log(2) \rfloor + 1
$$
iterations, our value of $j_m$ will be less than $j'$.  
So, after at most 
$$
T\ <\ 20 j^{1/2} \log j
$$ 
steps we reach our 
$$
\delta_T\ >\ 2/3\ {\rm or\ } j_T\ <\ j_0,
$$  
which finishes the verification of (\ref{inv5}), and so finishes the proof of our theorem.

\section{Proof of Proposition \ref{prop0}}

The proof of this proposition is fairly complex, and itself requires two
long lemmas, both of which are proved in seperate subsections within
this section.  The next subsection contains the statements of these lemmas.

\subsection{Preliminary Lemmata}

\begin{lemma} \label{prop1} 
We begin by supposing that $B > 1$ is some integer constant, that $n$ is
sufficiently large (as a function of $B$), and that 
$f : \F \to [0,1]$ satisfies 
$$
0\ <\ \theta\ :=\ \EE(f)\ <\ 1/4.
$$

If there exists an index 
$$
1\ \leq\ \ell\ \leq\ F/B
$$
such that 
$$
{\rm if\ }|\hat f(a_\ell)|\ =\ \gamma F,\ \ {\rm then\ \ } 
|\hat f(a_{B \ell})|\ <\ \theta \gamma^2 F,
$$
then one of the following two conclusions must hold:
\bigskip

$\bullet$ Either
\begin{equation} \label{prop1_conclusion}
\LL(f)\ >\ \theta \gamma^2/4;
\end{equation}
\bigskip

$\bullet$ or, if we let 
$$
R_1\ :=\ \{a_1,...,a_\ell\},\ {\rm and\ } R_2\ :=\ \{a_1,...,a_{B\ell}\},
$$
then there exists $t \in \F$ such that 
$$
|R_2 \cap (R_2 + t)|\ \geq\ (\ell/B)^{1/2}.
$$
\end{lemma}

\subsection{Body of the proof of Proposition \ref{prop0}} 

We begin by noting that we may assume that $\theta < 1/4$, since otherwise
Meshulam's theorem \cite{meshulam} 
implies that $S$ contains a three-term arithmetic progression
once $n$ is sufficiently large.

Next, as in the hypotheses of our Proposition, we assume that
that 
$$
\delta\ >\ \eps \ >\ 0, 
$$
and we let 
$$
B\ :=\ 50.
$$

The proof of our Proposition amounts to verifying that we can perform 
a certain type of ``Roth iteration'', where at each step the number of ``large'' 
Fourier coefficients decreases a lot, due to the additive structure of 
$R_2$ elucidated in Lemma \ref{prop1}. 

\subsubsection{Two cases}

Now, let
$$
k\ :=\ \lceil j/B \rceil,\ k\ >\ k_0.
$$
(which forces $j$ to be sufficiently large in terms of $\eps$), and then define
$\gamma > 0$ via the relation
$$
|\hat f(a_k)|\ =\ \gamma F.
$$
Then, we either have that 
\begin{equation} \label{tk}
|\hat f(a_{Bk})|\ <\ \gamma^3 F,
\end{equation}
or we don't.  
\bigskip

\noindent {\bf Case 1 (reverse inequality holds).}
If the reverse inequality holds, by which we mean that
$$
|\hat f(a_{B k})|\ \geq\ \gamma^3 F,
$$
then it follows that for $k$ sufficiently large,
$$
|\hat f(a_k)|\ \leq\ \theta^{j^{1/2 + \delta}/3} F\ <\ \theta^{ 2 k^{1/2 + \delta}} F,
$$
which is one of the conclusions of our Proposition.
\bigskip

\noindent {\bf Case 2 (inequality (\ref{tk}) holds).}
On the other hand, if (\ref{tk}) holds, then so long as $k > k_0$ and $n$ is sufficiently
large, we will have that the hypotheses of Lemma \ref{prop1} are met for 
$\ell = k$.  So, one or the other of the conclusions of Lemma \ref{prop1} must hold.

If the first conclusion of Lemma \ref{prop1} holds, 
then we have that $S$ contains a three-term arithmetic
progression since we have from (\ref{faj_bound}), along with the fact 
$k \leq j-1$, that
$$
\LL(f)\ >\ \theta \gamma^2/4\ >\ \theta (2 F^{-1/2})^2/4\ >\ \theta F^{-1}.
$$

If the second conclusion of Lemma \ref{prop1} holds, then we have that 
$$
R_2\ :=\ \{a_1,...,a_{B k}\}
$$
satisfies
$$
|R_2 \cap (R_2 + t)|\ \geq\ k^{1/2},
$$
for some $t \in \F$.  So, since 
$$
j\ \leq\ Bk\ <\ j + B,
$$
we deduce that if 
$$
R\ :=\ \{a_1,...,a_j\},
$$
then
\begin{equation} \label{Rt}
|R \cap (R+t)|\ \geq\ k^{1/2} - B\ \geq\ (j/B)^{1/2} - B.
\end{equation}

We now initiate another sub-subsection to expound upon this last case.

\subsubsection{Construction of the function $h$, and conclusion of the proof}

Let $t$ be as in (\ref{Rt}), and then define
$$
V\ :=\ \{ v \in \F\ :\ v \cdot t\ =\ 0\};
$$
that is, $V$ is the orthogonal complement of the one-dimensional subspace 
generated by $t$.   Next, suppose that $x$ is some multiple of $t$, and then
define the function
$$
g(n)\ :=\ f(n-x) V(n),
$$
where $V(n)$ is just the indicator function for $V$.  
Since $V$ is isomorphic as a vector space to $\F_{p^{n-1}}$, say the isomorphism is
$$
\varphi\ :\ \F_{p^{n-1}}\ \to\ V,
$$
then the new function
$$
h(n)\ =\ (g \circ \varphi)(n)
$$
satisfies
$$
h\ :\ \F_{p^{n-1}}\ \to\ [0,1]
$$
and if we let 
$$
T\ =\ {\rm support}(h),
$$
then $T$ has a non-trivial three-term arithmetic progression implies that $g$,
and therefore $f$, both do as well.

Passing from $f$ to this new function $h$ defined on a smaller dimensional 
space, constitutes one Roth-Meshulam interation.    We now consider
what the Fourier coefficients of $h$ look like:
As is well known (and easy to show), the Fourier coefficients of $h$ are
the same as those of $g$ (when the Fourier transform is restricted to $V$), since
vector space isomorphisms preserve Fourier spectra.   Thus, to compute the
largest Fourier coefficients of $h$, we just need to compute those of $g$.  
With a little work one can see that for $v \in V$,
\begin{equation} \label{gv}
\hat g(v)\ =\ p^{-1} \SIGMA_{u \in \F_p} e^{2\pi i u x/p} \hat f(v + tu);
\end{equation}
that is, $\hat g(v)$ is some sort of weighted average over $\hat f(v+tu)$
where $u$ ranges over $\F_p$.  

Let us assume that $x$ is chosen so that $\hat g(0)$ is maximal, and
therefore satisfies $\EE(h) \geq \theta$.  Such $x$
exists becasue the average of $\hat g(0)$ over all $x$ is $\theta p^{n-1}$.

Next, let us consider how many Fourier coefficients $\hat g(v)$ satisfy
$$
|\hat g(v)|\ >\ \theta^{j^{1/2 + \delta}} p^n.
$$
Clearly any such $v$ must have the property that at least one of 
$$
|f(v)|,\ |f(v+t)|,\ ...,\ {\rm or\ } |f(v + (p-1)t)|\ >\ \theta^{j^{1/2 + \delta}} p^n.
$$
Furthermore, since we are assuming that (\ref{Rt}) holds, we must have that 
for at least $(j/B)^{1/2} - B$ of these values $v \in V$, 
the sum in (\ref{gv}) contains at least two elements from
$$
R\ :=\ \{a_1,...,a_j\}.
$$
What that means is that there are a lot fewer $v$ where $|\hat g(v)|$ is large,
than there were places $a$ where $|\hat f(a)|$ is large.  In fact, we will have that 
if we write
$$
\F_{p^{n-1}}\ :=\ \{b_1,...,b_G\},\ {\rm where\  } G\ :=\ p^{n-1},
$$
where the $b_i$ are arranged so that
$$
|\hat h(b_1)|\ \geq\ |\hat h(b_2)|\ \geq\ \cdots\ \geq\ |\hat h(b_G)|,
$$
then for any
$$
\ell\ \geq\ j\ -\ (j/B)^{1/2}+ B,
$$
we will have that 
$$
|\hat h(b_{\ell})|\ \leq\ p \theta^{j^{1/2 + \delta}} G.
$$

This finishes the proof of our Proposition.

\subsection{Proof of Lemma \ref{prop1}}

Suppose that the hypotheses of the Lemma hold.
We will first establish the following claim.
\bigskip

\noindent {\bf Claim.}  Suppose that there exists a pair of points 
$$
b_1,\ b_2\ \in\ R_1,
$$
such that the only triple of the form
$$
-2a,\ a + b_1,\ a + b_2
$$
lying in $R_2$, is 
$$
0,\ b_1,\  b_2.
$$   
Then, we must have that (\ref{prop1_conclusion}) holds.
\bigskip

\noindent {\bf Proof of the claim.}  We begin with the basic fact that 
$$
\LL(f)\ =\ F^{-3} \SIGMA_a \hat f(a)^2 \hat f(-2a).
$$
What we will do is modify the product of the two 
$\hat f(a)$ in the $\hat f(a)^2$ to $\hat f_1(a) \hat f_2(a)$ in 
such a way that we can produce a lower bound for $\LL(f)$, while at the same 
time making use of the hypothesis of the claim.
\bigskip

These new functions are 
$$
f_1(m)\ :=\ f(m) e^{2\pi i b_1\cdot m/p},
$$
and
$$
f_2(m)\ :=\ f(m) e^{2\pi i b_2\cdot m/p}.
$$
Clearly, 
$$
{\rm support}(f_1)\ =\ {\rm support}(f_2)\ =\ {\rm support}(f);
$$
and so, we must have that 
$$
\LL(f)\ \geq\ F^{-3} \Bigl | \SIGMA_a \hat f_1(a) \hat f_2(a) \hat f(-2a) \Bigr |.
$$
\bigskip

Let us now look at this sum over $a$ a little more closely:  First, we observe that
$$
\hat f_1(a)\ =\ \SIGMA_m f(m) e^{2\pi i m\cdot (b_1 + a)/p}\ =\ \hat f(a+b_1),
$$
and
$$
\hat f_2(a)\ =\ \SIGMA_m f(m) e^{2\pi i m \cdot (b_2 + a)/p}\ =\ \hat f(a+b_2).
$$
Thus,
$$
\LL(f)\ \geq\ F^{-3} \Bigl | \SIGMA_a \hat f(a+b_1) \hat f(a+b_2) \hat f(-2a) \Bigr |.
$$

From this, and the hypotheses of our claim, we can easily deduce that
\begin{eqnarray*}
\LL(f)\ &\geq&\ F^{-3} |\hat f(b_1) \hat f(b_2) \hat f(0)|\ -\ 
3 F^{-3} \sup_{a \in \F \setminus R_2} |\hat f(a)| \SIGMA_a |\hat f(a)|^2 \\
&\geq&\ \theta \gamma^2\ -\ 3 ( \theta \gamma^2) \theta. 
\end{eqnarray*}

This then proves the claim, as $\theta < 1/4$.

$\hfill$ $\blacksquare$
\bigskip

From this claim we easily deduce that either (\ref{prop1_conclusion}) holds,
or else for every pair
$$
b_1, b_2\ \in\ R_1
$$
there exists $a$ such that 
$$
-2a, a + b_1, a + b_2\ \in\ R_2.
$$
So, by the pigeonhole principle, 
either (\ref{prop1_conclusion}) holds, or else there exists $a$ such that
$-2a \in R_2$ and  
\begin{eqnarray*}
|\{ (b_1,b_2) \in R_1 \times R_1\ :\ a + b_1, a+b_2 \in R_2\}|\ &\geq&\ \ell^2 (B\ell)^{-1} \\
&=&\ \ell/B.
\end{eqnarray*}
In other words,
$$
|\{ b \in R_1\ :\ a + b \in R_2\}|\ \geq\ (\ell/B)^{1/2}.
$$ 
So,
$$
|R_2 \cap (R_2 + a)|\ \geq\ |R_2 \cap (R_1 + a)|\ \geq\ (\ell/B)^{1/2},
$$
and the Lemma follows for $t = a$.

\section{Proof of Proposition \ref{prop30}}

First, we require the following:
\bigskip

\begin{lemma} \label{subspace_lemma}  For $1 \leq j < F^{1/8}$, 
there exists an additive subgroup $V$ of $\F_{p^n}$ having 
dimension $\lfloor 3n/4 \rfloor$, so that all the cosets
$$
a_1 + V, ..., a_j + V, -2a_1 + V,..., -2a_j + V,
$$
are distinct.
\end{lemma}
The proof of this lemma can be found in a later subsection within this section.
\bigskip

Now we define an auxiliary function $g$ to be 
$$
g(m)\ =\ (f W * V)(m)\ =\ \SIGMA_{b \in V} f(m-b) W(m-b),
$$
where $V(m)$ is the indicator function of the subsapce $V$ given in Lemma 
\ref{subspace_lemma}, and where $W(m)$ is the indicator function of 
$W = V^\perp$.  A simple calculation reveals that 
$$
\hat g(a)\ =\ \left \{ \begin{array}{rl}  \SIGMA_{v \in V} \hat f(a + v),\ & {\rm if\ } a\in W; \\
0,\ & {\rm if\ } a \notin W.
\end{array}\right.
$$

We now introduce some additional notation:  Given an $a \in F$, we write 
$a$ uniquely as 
$$
a\ =\ w(a) + v(a),\ {\rm where\ } w(a) \in W,\ {\rm and\ } v(a) \in V.
$$
From the conclusion of Lemma \ref{subspace_lemma} above, we have that 
for $i=1,...,j$, 
\begin{eqnarray} \label{fg_diff}
| \hat g(w(a_i))\ -\ \hat f(a_i)|\ &=&\ \left | \SIGMA_{v \in V \atop v \neq v(a_i)} 
\hat f(w(a_i) + v) \right | \nonumber \\
&\leq&\ |V| \sup_{v \in V \atop v \neq v(a_i)} |\hat f(w(a_i)+v)| \nonumber \\ 
&\leq&\ |V|\cdot |\hat f(a_{j+1})| \nonumber \\
&<&\ \theta^2 |V|/32.
\end{eqnarray}
Here, we have used the conclusion of Lemma \ref{subspace_lemma}, which
implies that  
$$
\{a_1,...,a_j\}\ \cap\ \{ w(a_i) + v\ :\ v \in V,\ v \neq v(a_i)\}\ =\ \emptyset.
$$
Also, note that one of the hypotheses of Proposition \ref{prop30} implies that
$$
|\hat f(a_{j+1})|\ \leq\ |\hat f(a_j)|\ <\ \theta^2/32.
$$

We likewise can deduce from Lemma \ref{subspace_lemma} that 
for $i=1,...,j$,
\begin{equation} \label{double_replace}
|\hat g(-2 w(a_i))\ -\ \hat f(-2a_i) |\ <\ \theta^2 |V|/32.
\end{equation}
On the other hand, we have that 
\begin{equation} \label{other_w}
{\rm if\ }w \in W,\ w \neq w(a_i)\ ({\rm for\ any\ }i=1,...,j),\ {\rm then\ } 
|\hat g(w)|\ <\ \theta^2|V|/32.
\end{equation}
\bigskip

Before pressing on, let us point out one conclusion of (\ref{fg_diff}) that we 
will use later on:
\begin{equation} \label{g_exp}
|\EE(g)\ -\ \EE(f)|\ =\ F^{-1} |\hat g(0) - \hat f(0)|\ <\ F^{-1} |V|\ =\ |W|^{-1}.
\end{equation}

Putting together the above observations we can deduce that 
\begin{eqnarray*}
\SIGMA_{m,d \in \F} g(m)g(m+d) g(m+2d)\ &=&\ 
F^{-1} \SIGMA_{w \in W} \hat g(w)^2 \hat g(-2w) \\ 
&=&\ F^{-1} \SIGMA_a \hat f(a)^2 \hat f(-2a)\ +\ E \\
&=&\ \SIGMA_{m,d \in \F} f(m)f(m+d)f(m+2d)\ +\ E,
\end{eqnarray*}
where $E$ is a certain error that can be computed through the use of 
the Cauchy-Schwarz inequality as follows.  

\subsection{The error $E$}

First, we observe that 
$$
F^{-1} \SIGMA_a \hat f(a)^2 \hat f(-2a)\ =\ 
F^{-1} \SIGMA_{i=1}^j \hat f(a_i)^2 \hat f(-2a_i)\ +\ E_1,
$$
where by Parseval and Cauchy-Schwarz,
$$
|E_1|\ \leq\ F^{-1} |\hat f(a_{j+1})| \SIGMA_a |\hat f(a)|^2\ <\ \theta^3 F/32.
$$
Next, we have that 
$$
F^{-1} \SIGMA_{i=1}^j \hat f(a_i)^2 \hat f(-2a_i)\ =\ 
F^{-1} \SIGMA_{i=1}^j \hat g(w(a_i)) \hat f(a_i) \hat f(-2a_i)\ +\ E_2,
$$
where by Parseval and Cauchy-Schwarz
$$
|E_2|\ =\ \left | F^{-1} \SIGMA_{i=1}^j (\hat f(a_i) - \hat g(w(a_i)) \hat f(a_i) \hat f(-2a_i) \right |
\ <\ \theta^3 |V| F/32. 
$$

Next, we replace another of the factors $\hat f(a_i)$ with $\hat g(a_i)$, incurring a small
error:   We have that 
$$
F^{-1} \SIGMA_{i=1}^j \hat g(w(a_i)) \hat f(a_i) \hat f(-2a_i)\ =\ 
F^{-1} \SIGMA_{i=1}^j \hat g(w(a_i))^2 \hat f(-2a_i)\ +| E_3|,
$$
where by Cauchy-Schwarz, Parseval, the fact that $0 \leq g(m) \leq 1$ and 
(\ref{g_exp}), we have that 
\begin{eqnarray*}
|E_3|\ =\ \left | F^{-1} \SIGMA_{i=1}^j \hat g(w(a_i)) (\hat f(a_i) - \hat g(w(a_i))) 
\hat f(-2a_i) \right |\ &<&\ \theta ( \theta + |W|^{-1}) |V| F \\
&<&\ \theta^3 |V| F/16.
\end{eqnarray*}

Next, we replace the $\hat f(-2a_i)$ with $\hat g(-2w(a_i))$ using 
(\ref{double_replace}), by first writing
$$
F^{-1} \SIGMA_{i=1}^j \hat g(w(a_i))^2 \hat f(-2a_i)\ =\ 
F^{-1} \SIGMA_{i=1}^j \hat g(w(a_i))^2 \hat g(-2w(a_i))\ +\ E_4,
$$
where
\begin{eqnarray*}
|E_4|\ &=&\ \left | F^{-1} \SIGMA_{i=1}^j \hat g(w(a_i))^2 (\hat f(-2a_i) - \hat g(-2w(a_i))) \right |
\\ 
&<&\ \theta^2(\theta + |W|^{-1})^2 |V| F/32 \\  &<&\ \theta^3 |V| F/16. 
\end{eqnarray*}

Finally, we consider the complete sum
$$
F^{-1} \SIGMA_{w \in W} \hat g(w)^2 \hat g(-2w)\ =\ F^{-1} \SIGMA_{i=1}^j 
\hat g(w(a_i))^2 \hat g(-2w(a_i))\ +\ E_5,
$$
where by (\ref{other_w}), Cauchy-Schwarz, and Parseval, we have that
$$
|E_5|\ <\ F^{-1} |V| (\theta + |W|^{-1}) \theta^2 F^2/32\ <\ \theta^3 |V| F/16.
$$
\bigskip

Combining the errors $E_1,E_2,E_3,E_4$ and $E_5$, we deduce that 
$$
|E|\ <\ \theta^3 |V| F/4.
$$

\subsection{Resumption of the proof of Proposition \ref{prop30}}

To finish the proof of the proposition, we derive a lower bound for $\LL(g)$, and
then a lower bound for $\LL(f)$:  First, observe that since $g$ is translation-invariant
by elements of $V$, we have that 
$$
\LL(g)\ \geq\ F^{-2} |V|^2 \SIGMA_{w \in W} g(w)^3\ \geq\ |W|^{-2} |W| (\theta - |W|^{-1})^3
\ >\ (2|W|)^{-1} \theta^3.  
$$
On the other hand, our bound on $E$ above guarantees that 
\begin{eqnarray*}
\LL(f)\ &\geq&\ \LL(g)\ -\ F^{-2} |E|\ \geq\ (2 |W|)^{-1} \theta^3\ -\ (4|W|)^{-1} \theta^3 \\
&=&\ (4 |W|)^{-1} \theta^3,
\end{eqnarray*}
which exceeds the trivial lower bound of $\theta |F|^{-1}$ so long as
$$
\theta^2\ >\ 4 |V|^{-1}. 
$$
In other words, our theorem holds, so long as 
$$
\theta\ >\ 2 p F^{-1/8}.
$$

\subsection{Proof of Lemma \ref{subspace_lemma}}

By basic properties of subspaces, 
the claim of the lemma is easily seen to be implied by the statement
\begin{equation} \label{intersection}
B \cap V\ =\ \emptyset,
\end{equation}
where 
$$
B\ :=\ \{a_{i_1} - a_{i_2}\ :\ 1 \leq i_1 < i_2 \leq j\}\ \cup\ \{2a_{i_1} + 
a_{i_2}\ :\ 1 \leq i_1 < i_2 \leq j\}.
$$
Note that 
$$
|B|\ \leq\ j(j-1).
$$
\bigskip

The idea will be to show that with positive probability, 
a randomly chosen subsapce $V$ of dimension 
$$
n'\ =\ \lfloor 3n/4 \rfloor,
$$
satisfies (\ref{intersection}).  To this end, we first observe that 
the probability that any non-zero element of $\F$ happens to lie in our 
ranomly chosen $V$ is the same as any other non-zero element of $\F$; so,
for any $x \in \F\setminus \{0\}$ we have that 
$$
{\rm Prob}(x \in V)\ =\ {|V| - 1 \over F - 1}.
$$
It follows that the probability that none of the element of $B$ happen to lie in $V$
is at least 
$$
1\ -\ j(j-1) {|V|-1 \over F - 1}\ >\ 1\ -\ { (F^{1/8})(F^{1/8} -1)(F^{3/4} - 1) \over F- 1}\ >\ 0.
$$
This completes the proof of the lemma.

\section{Proof of Proposition \ref{prop3}}

First, we may assume that $\theta < 1/p$, since if $\theta \geq 1/p$ we have
by Meshulam's theorem \cite{meshulam} that for $n$ sufficiently large 
that the set $S$ contains a three-term arithmetic progression.

The proof of the proposition will be very different according as to whether 
$$
j\ \leq\ 3n/2,\ {\rm or\ } 3n/2\ \leq\ j\ <\ F^{1/8}.
$$

\subsection{Case 1:  $j \leq 3n/2$}

As the title of this subsection suggests we will assume that 
$$
j\ \leq\ 3n/2.
$$

Let 
$$
V\ :=\ {\rm span}(a_1,...,a_j).
$$
Using the fact that $a_2 = -a_3$, $a_4 = -a_5$, and so on, we deduce that
$$
V\ =\ \left \{ \begin{array}{rl} {\rm span}(a_2,a_4,...,a_j),\ & {\rm if\ } j\ {\rm even}; \\
{\rm span}(a_2,a_4,...,a_{j-1}),\ & {\rm if\ } j\ {\rm odd}. \end{array}\right.
$$
Note that in either case, we will have that 
$$
{\rm dim}_{\F_p}(V)\ \leq\ 3n/4.
$$
Next, define $W$ to be the orthogonal complement of $V$, and let
$$
g(n)\ =\ f(n-x)W(n),
$$
where $W(n)$ is the indicator function for $W$.  By averaging we can clearly choose
$x$ so that 
$$
\EE(g)\ \geq\ |V|^{-1} \theta;\ {\rm or,\ equivalently,\ } 
\hat g(0)\ \geq\ \theta |W|.
$$
Now, if we let $T$ be the support of $g$, then $T$ contains a non-trivial 
three-term arithmetic progression clearly implies that $S$ does as well; and,
to decide whether $T$ has three-term progressions, we compute 
\begin{eqnarray} \label{3apg}
\SIGMA_{a,d \in W} g(a)g(a+d)g(a+2d)\ &=&\ |W|^{-1} \SIGMA_{b \in W} 
\hat g(b)^2 \hat g(-2b) \nonumber \\
&=&\ |W|^{-1} \hat g(0)^3 -\ E,
\end{eqnarray}
where the error $E$ satisfies
$$
|E|\ \leq\ M |W|^{-1} \SIGMA_{b \in W} |\hat g(b)|^2\ \leq\ 
M  \hat g(0),
$$
where
$$
M\ :=\ \sup_{b \in W \setminus \{0\}} |\hat g(b)|.
$$
So, the quantity in (\ref{3apg}) is at least
\begin{equation} \label{3apg2}
\hat g(0) \left ( |W|^{-1} \hat g(0)^2 - M \right ).
\end{equation}
To bound $M$ from above, we will need a formula for the Fourier transfrom
$\hat g(b)$, and such a formula (which is easy to show) is
$$
\hat g(b)\ =\ |V|^{-1}\SIGMA_{v \in V} e^{2\pi i x\cdot v/p} \hat f(b + v).
$$
Now, if $b = 0$, then this sum will include all the $\hat f(a_i)$ for
$i=1,...,j$; and, if $b \in W \setminus \{0\}$, then the sum includes none
of these numbers, which implies that
$$
{\rm for\ } b \in W \setminus \{0\},\ |\hat g(b)|\ \leq\ |\hat f(a_{j+1})|\ \leq\ 
\theta^{j+2} F/2\ <\ \theta^2 p^{-j} F/2.
$$
Thus, 
$$
M\ <\ p^{-j} \theta^2 F/2\ <\ \theta^2 W/2,
$$
and it follows that the quantity in (\ref{3apg2}) is at least
$$
\hat g(0) \left ( |W|^{-1} \hat g(0)^2 - \theta^2|W|/2 \right )\ >\ 
\theta^3 |W|^2/2.
$$
In order for $S$ to contain a non-trivial three-term arithmetic progression, we need 
that this last quantity exceeds $\theta |W|$, and it does provided that
$$
\theta\ >\ 2 |W|^{-1/2}\ =\ 2 (p^{n - {\rm dim}(V)} )^{-1/2}\ >\ 2 p^{-n/8}\ =\ 2 F^{-1/8}.
$$

\subsection{Case 2 : $3n/2\ <\ j\ <\ F^{1/8}$}

To handle this case we will apply Proposition \ref{prop30}:  First, from the hypothesis
of our Proposition \ref{prop3}, along with the assumption $j > 3n/2$ we have for 
$n$ sufficiently large, and $\theta < 1/p$, that 
\begin{equation} \label{ajc}
|\hat f(a_j)|\ <\ \theta^{j+2} F/2\ <\ \theta^2/32.
\end{equation}
Of course if $\theta \geq 1/p$, then we know from Meshulam's theorem \cite{meshulam}
that for $n$ large enough, $S$ contains a three-term arithmetic progression.

So, since (\ref{ajc}) holds we have by Proposition \ref{prop30} that $S$ contains 
a three-term arithmetic progression, and our proposition is proved.


\begin{thebibliography}{999}

\bibitem{croot} E. Croot, {\it On the decay of the Fourier transform and three term
arithmetic progressions}, Online Jour. of Analy. Comb. {\bf 2} (2007).

\bibitem{meshulam} R. Meshulam, {\it On subsets of finite abelian groups with no
3-term arithmetic progressions}, J. Comb. Theory Ser. A {\bf 71} (1995), 168-172.

\bibitem{shkredov} I. D. Shkredov, {\it On sets of large exponential sum}, preprint on 
the arxives. 

\end{thebibliography}
\end{document}